\documentclass[11pt, reqno]{amsart}
\usepackage{amssymb}
\usepackage{amscd}
\usepackage{verbatim,ifthen}
\usepackage{color}
\usepackage{latexsym}
\usepackage{tikz}
\usepackage{tikz-cd}
\usepackage{mathrsfs}
\usepackage{wrapfig}
\usetikzlibrary{shapes}
\usepackage{color}
\usetikzlibrary{arrows.meta}
\usepackage{bbm}
\usetikzlibrary{matrix}
\usetikzlibrary{calc}
\usetikzlibrary{arrows,intersections}
\usepackage{pgfplots}
\usepackage{multicol}
\usepackage{array}
\newcolumntype{M}[1]{>{\centering\arraybackslash}m{#1}}

\usepackage[colorlinks, linkcolor=black, citecolor=magenta, linktocpage]{hyperref}
\addtocontents{toc}{\setcounter{tocdepth}{1}}

\usepackage{amsmath}

\numberwithin{equation}{section}

\let\oldtocsection=\tocsection
 
\let\oldtocsubsection=\tocsubsection

\renewcommand{\tocsection}[2]{\hspace{0em}\oldtocsection{#1}{#2}}
\renewcommand{\tocsubsection}[2]{\hspace{1em}\oldtocsubsection{#1}{#2}}

\def\Xint#1{\mathchoice
{\XXint\displaystyle\textstyle{#1}}%
{\XXint\textstyle\scriptstyle{#1}}%
{\XXint\scriptstyle\scriptscriptstyle{#1}}%
{\XXint\scriptscriptstyle\scriptscriptstyle{#1}}%
\!\int}
\def\XXint#1#2#3{{\setbox0=\hbox{$#1{#2#3}{\int}$ }
\vcenter{\hbox{$#2#3$ }}\kern-.6\wd0}}

\def\dashint{\Xint-}
\usepackage{eucal}
\usepackage{calc}  
\usepackage{enumitem} 
\usepackage{tensor}
\usepackage{graphicx,wrapfig,lipsum}
\usepackage{etoolbox}
\usepackage{marginnote}
\usepackage{lipsum}
\makeatletter
\patchcmd{\@mn@margintest}{\@tempswafalse}{\@tempswatrue}{}{}
\patchcmd{\@mn@margintest}{\@tempswafalse}{\@tempswatrue}{}{}
\reversemarginpar 
\makeatother
\usepackage{scrextend}

\makeatletter
\DeclareRobustCommand\widecheck[1]{{\mathpalette\@widecheck{#1}}}
\def\@widecheck#1#2{%
    \setbox\z@\hbox{\m@th$#1#2$}%
    \setbox\tw@\hbox{\m@th$#1%
       \widehat{%
          \vrule\@width\z@\@height\ht\z@
          \vrule\@height\z@\@width\wd\z@}$}%
    \dp\tw@-\ht\z@
    \@tempdima\ht\z@ \advance\@tempdima2\ht\tw@ \divide\@tempdima\thr@@
    \setbox\tw@\hbox{%
       \raise\@tempdima\hbox{\scalebox{1}[-1]{\lower\@tempdima\box
\tw@}}}%
    {\ooalign{\box\tw@ \cr \box\z@}}}
\makeatother
\title{On the Weighted Orthogonal Ricci curvature}
\author{Kyle Broder and Kai Tang}

\begin{document}

\maketitle

\begin{abstract}
We introduce the weighted orthogonal Ricci curvature -- a two-parameter version of Ni--Zheng's orthogonal Ricci curvature. This curvature serves as a very natural object in the study of the relationship between the Ricci curvature(s) and the holomorphic sectional curvature. In particular, in determining optimal curvature constraints for a compact K\"ahler manifold to be projective. In this direction, we prove a number of vanishing theorems using the weighted orthogonal Ricci curvature(s) in both the K\"ahler and Hermitian category.
\end{abstract}

\section{Introduction}

Let $(M^n, \omega)$ be a compact K\"ahler manifold. One of the central problems in modern complex geometry is the relationship between the Ricci curvature $\text{Ric}_{\omega}$ and the holomorphic sectional curvature $\text{HSC}_{\omega}$. The Wu--Yau theorem \cite{WuYau1, WuYau2, TosattiYang, DiverioTrapani, YangZhengRBC, Nomura, Tang0, LeeStreets, BroderSBC, BroderSBC2} is one such manifestation of this relationship: If $M$ supports a metric of negative holomorphic sectional curvature, then $M$ supports a (possibly different) metric with negative Ricci curvature. It is known from Hitchin's metrics on Hirzebruch surfaces \cite{Hitchin}, however, that $\text{HSC}_{\omega} >0$ does not imply $\text{Ric}_{\omega} >0$.

In an attempt to further understand this relationship, Ni--Zheng \cite{NiZhengComparison} have studied the orthogonal Ricci curvature (which first appeared under the name ``anti-holomorphic Ricci curvature" in \cite{MiquelPalmer}): \begin{eqnarray*}
\text{Ric}_{\omega}^{\perp} : T^{1,0} M \to \mathbb{R}, \hspace{1cm} \text{Ric}_{\omega}^{\perp}(X) \ : = \ \frac{1}{| X |_{\omega}^2}\text{Ric}_{\omega}(X, \overline{X}) - \text{HSC}_{\omega}(X).
\end{eqnarray*}

It was observed in \cite{NiZhengComparison, HuangTam, Niu} that $\text{Ric}_{\omega}^{\perp}$ is, in an appropriate sense, the trace of the orthogonal bisectional curvature $\text{HBC}_{\omega}^{\perp}$ (i.e., the restriction of the bisectional curvature to pairs of orthogonal vectors). This justifies the terminology orthogonal Ricci curvature (c.f., \nameref{Remark2.7}). 

It is suspected (c.f., \cite{NiZhengComparison}) that the orthogonal Ricci curvature should give a sharp measure of  the projectivity of compact K\"ahler manifolds. For instance, $\text{Ric}_{\omega} >0$ implies that the anti-canonical bundle $K_M^{-1}$ is ample, hence $M$ is projective. Recently, it was shown by Yang \cite{YangYauConjecture} that $\text{HSC}_{\omega} >0$ implies $h^{p,0}=0$ for all $1 \leq p \leq n$. Ni \cite{NiRational} showed that any compact K\"ahler manifold with $\text{Ric}_{\omega}^{\perp}>0$ must be projective and simply connected (also see \cite{NiZhengComparison}).

With this aim of generating a sharp curvature constraint for compact K\"ahler manifolds to be projective, we introduce the following \textit{weighted orthogonal Ricci curvature}:

\subsection*{Definition 1.1}
Let $(M^n, \omega)$ be a compact K\"ahler manifold. For $\alpha, \beta \in \mathbb{R}$, we define the \textit{weighted orthogonal Ricci curvature} to be the function $$\text{Ric}_{\alpha, \beta}^{\perp} : T^{1,0}M \to \mathbb{R}, \hspace*{1cm} \text{Ric}_{\alpha,\beta}^{\perp} : = \frac{\alpha }{| X |_{\omega}^2} \text{Ric}_{\omega}(X, \overline{X}) + \beta\text{HSC}_{\omega}(X).$$ 

\subsection*{Remark 1.2}
Of course, when $\alpha = 1$, $\beta=-1$, we recover the orthogonal Ricci curvature $\text{Ric}^{\perp}$. Chu--Lee--Tam \cite{ChuLeeTam} showed that, for $\alpha>0$ and $\beta>0$, a compact K\"ahler manifold with $\text{Ric}_{\alpha, \beta}^{\perp}<0$ must be projective with ample canonical bundle, extending the aforementioned Wu--Yau theorem. Further, they also show that a compact K\"ahler manifold with $\text{Ric}_{\alpha, \beta}^{\perp}>0$ for $\alpha>0$, $\beta>0$, is projective and simply connected.

By considering the weighted orthogonal Ricci curvature, in place of merely the orthogonal Ricci curvature, natural questions concerning the \textit{extent} to which the holomorphic sectional curvature and Ricci curvature are related easily proliferate. For instance, the following question is of tremendous interest: 

\subsection*{Question 1.3}
Determine all $(\alpha, \beta) \in \mathbb{R}^2$ such that a compact K\"ahler manifold with $\text{Ric}_{\alpha, \beta}^{\perp} >0$ or $\text{Ric}_{\alpha,\beta}^{\perp}<0$ is projective. \\

In this direction, we have the following theorem:

\subsection*{Theorem 1.4}
Let $(M^n, \omega)$ be a compact K\"ahler manifold with $\text{Ric}_{\alpha, \beta}^{\perp} >0$ for some $\alpha > 0 > \beta$. If, moreover, $3 \alpha + 2 \beta >0$, then $M$ is projective.  \\

In a similar spirit, we have the following natural question: 

\subsection*{Question 1.5}
Determine all $(\alpha, \beta) \in \mathbb{R}^2$ such that a compact K\"ahler manifold with $\text{Ric}_{\alpha, \beta}^{\perp} >0$ satisfies $h^{p,0}=0$. \\

We have the following partial answer to the above question:

\subsection*{Theorem 1.6}
Let $(M^n, \omega)$ be a compact K\"ahler manifold with $\text{Ric}_{\alpha, \beta}^{\perp} >0$ for some $\alpha >0$ and $\beta <0$. If, moreover, $(p+1) \alpha + 2 \beta >0$, then $h^{p,0}=0$ for all $1 \leq p \leq n$. In particular, $M$ is projective.\\

We will exhibit a number of theorems of this type in this manuscript. The utility of such theorems is also seen if one has knowledge of the Hodge numbers. Indeed, we can show that certain manifolds do not support metrics with particular relations on the Ricci and holomorphic sectional curvature. To state one such instance of this, let us first observe that we can define Hermitian extensions of the weighted orthogonal Ricci curvature:

\subsection*{Definition 1.7}
Let $(M^n, \omega)$ be a Hermitian manifold. For $1 \leq k \leq 4$, we define the \textit{weighted $k$th orthogonal Ricci curvature} of $\omega$ to be the function \begin{eqnarray*}
\text{Ric}_{\alpha, \beta}^{(k)} : T^{1,0}M \to \mathbb{R}, \hspace*{1cm} \text{Ric}_{\alpha, \beta}^{(k)}(X) \ : = \  \frac{\alpha}{| X |_{\omega}^2} \text{Ric}_{\omega}^{(k)}(X, \overline{X}) + \beta \text{HSC}_{\omega}(X).
\end{eqnarray*}

\hfill

Since we know the Hodge numbers of the Iwasawa threefold, we can prove the following:

\subsection*{Theorem 1.8}
There is no balanced metric on the Iwasawa threefold such that $\text{Ric}_{\alpha, \beta}^{(k)} > 0$ for $\alpha > 0 > \beta$ and $3\alpha + 2 \beta >0$.

\section{Some Reminders of Curvature in Complex Geometry}
Let $(M^n, \omega)$ be a Hermitian manifold. The Chern connection $\nabla$ on $T^{1,0} M$ is the unique Hermitian connection whose torsion has vanishing $(1,1)$--part. Fix a point $p \in M$. In a local coordinate frame $\left \{ \frac{\partial}{\partial z_i} \right \}$ of $T_p^{1,0}M$, the components of the Chern curvature tensor $R$ read \begin{eqnarray*}
R_{i \overline{j} k \overline{\ell}} &=& - \frac{\partial^2 g_{k \overline{\ell}}}{\partial z_i \partial \overline{z}_j} + g^{p \overline{q}} \frac{\partial g_{k \overline{q}}}{\partial z_i} \frac{\partial g_{p \overline{\ell}}}{\partial \overline{z}_j}.
\end{eqnarray*}

\subsection*{Reminder 2.1: Chern--Ricci and scalar curvatures}
The Chern curvature is an $\text{End}(T^{1,0}M)$--valued $(1,1)$--form. The first Chern--Ricci curvature is the contraction over the endomorphism part: \begin{eqnarray*}
\text{Ric}^{(1)}_{\omega} &=& \sqrt{-1} \text{Ric}^{(1)}_{i \overline{j}} dz_i \wedge d\overline{z}_j \ = \ \sqrt{-1} g^{k \overline{\ell}} R_{i \overline{j} k \overline{\ell}} dz_i \wedge d\overline{z}_j,
\end{eqnarray*}

and is a $(1,1)$--form representing the first Chern class the anti-canonical bundle $c_1(K_M^{-1})$. The second Chern--Ricci curvature is a contraction over the $(1,1)$--part: \begin{eqnarray*}
\text{Ric}^{(2)}_{\omega} & =& \text{Ric}_{k \overline{\ell}}^{(2)} \ = \ g^{i \overline{j}} R_{i \overline{j} k \overline{\ell}}.
\end{eqnarray*}

Similarly, the third and fourth Chern--Ricci curvature are defined: $$\text{Ric}_{\omega}^{(3)} \ = \ \text{Ric}_{k \overline{j}}^{(3)} \ = \ g^{i \overline{\ell}} R_{i \overline{j} k \overline{\ell}}, \hspace*{1cm} \text{Ric}_{\omega}^{(4)} \ = \ \text{Ric}_{i \overline{\ell}}^{(4)} \ = \ g^{k \overline{j}} R_{i \overline{j} k \overline{\ell}}.$$ 

The contraction $$\text{Scal}_{\omega} \ : = \ g^{i \overline{j}} \text{Ric}_{i \overline{j}}^{(1)} \ = \ g^{i \overline{j}} g^{k \overline{\ell}} R_{i \overline{j} k \overline{\ell}}$$ is the Chern scalar curvature. The contraction $$\widehat{\text{Scal}}_{\omega} \ : = \ g^{i \overline{\ell}} \text{Ric}_{k \overline{j}}^{(3)} \ = \ g^{i \overline{\ell}} g^{k \overline{j}} R_{i \overline{j} k \overline{\ell}}$$ will be referred to as the alterred Chern scalar curvature.

\subsection*{Remark 2.2}
If the Hermitian metric is K\"ahler, then the Ricci curvatures all coincide. Of course, this implies that the scalar curvatures coincide if the metric is K\"ahler, too. In fact, these statements are true in the more general setting of K\"ahler-like metrics. That is, a metric is said to be K\"ahler-like \cite{YangZhengCurvature} if the Chern curvature tensor satisfies the symmetries of the K\"ahler curvature tensor. The Iwasawa threefold shows that a K\"ahler-like manifold is not necessarily K\"ahler. \\

For convenience, let us restate the definition of the curvatures of primary interest in this manuscript: 

\subsection*{Definition 2.3}
Let $(M^n, \omega)$ be a compact Hermitian manifold. For $1 \leq k \leq 4$, we define the \textit{weighted $k$th orthogonal Ricci curvature} of $\omega$ to be the function \begin{eqnarray*}
\text{Ric}_{\alpha, \beta}^{(k)} : T^{1,0}M \to \mathbb{R}, \hspace*{1cm} \text{Ric}_{\alpha, \beta}^{(k)}(X) \ : = \  \frac{\alpha}{| X |_{\omega}^2} \text{Ric}_{\omega}^{(k)}(X, \overline{X}) + \beta \text{HSC}_{\omega}(X).
\end{eqnarray*}

In particular, if the metric is K\"ahler\footnote{Or more generally, K\"ahler-like.}, we can speak of the weighted orthogonal Ricci curvature:

\subsection*{Definition 2.4}
Let $(M^n, \omega)$ be a compact K\"ahler manifold. For $\alpha, \beta \in \mathbb{R}$, we define the \textit{weighted orthogonal Ricci curvature} to be the function $$\text{Ric}_{\alpha, \beta}^{\perp} : T^{1,0}M \to \mathbb{R}, \hspace*{1cm} \text{Ric}_{\alpha,\beta}^{\perp} : = \frac{\alpha }{| X |_{\omega}^2} \text{Ric}_{\omega}(X, \overline{X}) + \beta \text{HSC}_{\omega}(X).$$ 

We say that $\text{Ric}_{\alpha, \beta}^{\perp}$ (or $\text{Ric}_{\alpha, \beta}^{(k)}>0$ for some $1 \leq k \leq 4$) is positive (negative, zero, quasi-positive, etc.) if this is true for some $(\alpha, \beta) \in \mathbb{R}^2$.\\

To justify the terminology, let us recall the Quadratic Orthogonal Bisectional Curvature introduced in \cite{WuYauZheng}:

\subsection*{Definition 2.5}
Let $(M^n, \omega)$ be a Hermitian manifold. The \textit{Quadratic Orthogonal Bisectional Curvature} (from now on, QOBC) is the function $$\text{QOBC}_{\omega} : \mathcal{F}_M \times \mathbb{R}^n \backslash \{ 0 \} \to \mathbb{R},  \hspace*{1cm} \text{QOBC}_{\omega}: (\vartheta, v) \mapsto \frac{1}{| v |_{\omega}^2} \sum_{\alpha, \gamma =1}^n R_{\alpha \overline{\alpha} \gamma \overline{\gamma}} (v_{\alpha} - v_{\gamma})^2,$$ where $R_{\alpha \overline{\alpha} \gamma \overline{\gamma}}$ denote the components of the Chern connection of $\omega$ with respect to the unitary frame $\vartheta$ (a section of the unitary frame bundle $\mathcal{F}_M$). \\

This curvature first appeared implicitly in \cite{BishopGoldberg} and is the Weitzenb\"ock curvature operator (in the sense of \cite{PetersenBook, PetersenWink, PetersenWinkHodge}) acting on real $(1,1)$--forms. See \cite{BroderGraph, BroderLinearAlgebra, BroderQOBC} for alternative descriptions of the QOBC. From \cite{LiWuZheng}, the QOBC is strictly weaker than the orthogonal bisectional curvature $\text{HBC}_{\omega}^{\perp}$ (the restriction of the holomorphic bisectional curvature $\text{HBC}_{\omega}$ to pairs of orthogonal $(1,0)$--tangent vectors). From \cite{GuZhang}, the K\"ahler--Ricci flow on a compact K\"ahler manifold with $\text{HBC}_{\omega}^{\perp} \geq 0$ converges to a K\"ahler metric $\text{HBC}_{\omega} \geq 0$. Hence, Mok's extension \cite{Mok} of the solution of the Frankel conjecture \cite{Mori, SiuYau} shows that all compact K\"ahler manifolds with $\text{HBC}_{\omega}^{\perp} \geq 0$ are biholomorphic to a product of Hermitian symmetric spaces (of rank $\geq 2$) and projective spaces. In particular, although $\text{HBC}_{\omega}^{\perp}$ is an algebraically weaker curvature notion than $\text{HBC}_{\omega}$, the positivity of $\text{HBC}_{\omega}^{\perp}$ does not generate new examples.

\subsection*{Remark 2.6}\label{Remark2.7}
Recall that the orthogonal Ricci curvature $\text{Ric}_{\omega}^{\perp}$, which first appeared in \cite{MiquelPalmer} under the name \textit{anti-holomorphic Ricci curvature}, is defined \begin{eqnarray*}
\text{Ric}_{\omega}^{\perp} : T^{1,0} M \to \mathbb{R}, \hspace{1cm} \text{Ric}_{\omega}^{\perp}(X) \ : = \ \frac{1}{| X |_{\omega}^2} \text{Ric}_{\omega}(X, \overline{X}) -  \text{HSC}_{\omega}(X).
\end{eqnarray*}

From \cite{NiZhengComparison, HuangTam,Niu}, we observe that, in a fixed unitary frame, we have \begin{eqnarray*}
\text{QOBC}_{\omega} &=& \sum_{i,k} R_{i \overline{i}} \rho_{k \overline{i}} \rho_{i \overline{k}} - \sum_{i,j} R_{i \overline{i} j \overline{j}} \rho_{i \overline{i}} \rho_{j \overline{j}} - \sum_{i \neq \ell \ \text{or} \ k \neq j} R_{i \overline{\ell} k \overline{j}} \rho_{\ell \overline{i}} \rho_{j \overline{k}}.
\end{eqnarray*}

Choosing $\rho_{i \overline{j}} =0$ for all $i,j$, except $\rho_{1 \overline{1}}$, we see that \begin{eqnarray*}
\text{QOBC}_{\omega} &=& (R_{1 \overline{1}} - R_{1 \overline{1} 1 \overline{1}}) | \rho_{1 \overline{1}} |^2 \ \geq \ 0.
\end{eqnarray*}

This implies that $\text{Ric}^{\perp} \geq 0$. Since $\text{Ric}_{\alpha, \beta}^{\perp} >0$ is defined \textit{for some} $\alpha, \beta \in \mathbb{R}$, it is immediate that $\text{Ric}^{\perp} > 0$ implies $\text{Ric}_{\alpha, \beta}^{\perp} >0$. Hence, the quadratic orthogonal bisectional curvature dominates the weighted orthogonal Ricci curvature $\text{Ric}_{\alpha, \beta}^{\perp}$.\\

In light of this remark, we pose the following:

\subsection*{Question 2.7}\label{Question2.8}
Determine all values of $(\alpha, \beta) \in \mathbb{R}^2$ such that $\text{QOBC}_{\omega} \geq 0$ implies $\text{Ric}_{\alpha, \beta}^{\perp} \geq 0$. \\

Let $\dashint_{\mathbb{S}^{2n-1}} : = \frac{1}{\text{vol}(\mathbb{S}^{2n-1})} \int_{\mathbb{S}^{2n-1}}$. If $(M^n, \omega)$ is a compact K\"ahler manifold, then we have \begin{eqnarray*}
\text{Scal}_{\omega} \ = \ 2n \phantom{.} \dashint_{\mathbb{S}^{2n-1}} \text{Ric}_{\omega}, \hspace*{0.7cm} \text{and} \hspace{0.7cm} \text{Scal}_{\omega} \ = \ n(n+1) \phantom{.} \dashint_{\mathbb{S}^{2n-1}} \text{HSC}_{\omega}.
\end{eqnarray*}

Hence, \begin{eqnarray*}
\dashint_{\mathbb{S}^{2n-1}} \text{Ric}_{\alpha,\beta}^{\perp} &=&  \frac{\alpha(n+1) + 2\beta}{2n(n+1)} \text{Scal}_{\omega},
\end{eqnarray*}

and subsequently, $\text{Ric}_{\alpha, \beta}^{\perp} >0$ implies $\text{Scal}_{\omega} >0$ if $\alpha(n+1) + 2 \beta >0$. Since $\text{QOBC}_{\omega} \geq 0$ implies $\text{Scal}_{\omega} \geq 0$ (see \cite{ChauTamHarmonic}), this argument gives a partial answer to \nameref{Question2.8}:

\subsection*{Proposition 2.8}
Let $(M^n, \omega)$ be a compact K\"ahler manifold with $\text{QOBC}_{\omega} \geq 0$. If $\text{Ric}_{\alpha, \beta}^{\perp} \geq 0$, then $\alpha(n+1) + 2 \beta \geq 0$. 

\section{Vanishing theorems}
In this section, we will prove a number of vanishing theorems for Hodge numbers, using the weighted orthogonal Ricci curvature(s).

An important technique that will be used throughout is Ni's technique of viscosity considerations (see \cite[$\S 3$]{NiRational}). Let us briefly describe this technique for the purposes needed here. Let $V$ a vector space, and let $\eta$ be a $k$--covector. We say that $\eta$ is \textit{simple} if there exist $v_1, ..., v_k$ such that $\eta = v_1 \wedge \cdots \wedge v_k$. Let $\mathcal{S}(1) : = \{ \eta : \eta \ \text{is simple with} \ \| \eta \| =1 \}$. Here, $\| \cdot \|$ is the norm induced by the scalar product defined on simple covectors $\eta = v_1 \wedge \cdots \wedge v_k$ and $\omega = u_1 \wedge \cdots \wedge u_k$ by the formula $$\langle \eta, \omega \rangle \ : = \ \det( \langle v_i, u_j \rangle)$$ and subsequently extended bilinearly to all $r$--covectors which are linear combinations of simple covectors. Define the \textit{comass} of a $k$--covector $\rho$: \begin{eqnarray*}
\| \rho \|_0 & : = & \sup_{\eta \in \mathcal{S}(1)} | \rho(\eta) |.
\end{eqnarray*}

Let $\sigma \in \Omega^{p,0}(M)$ be a holomorphic $(p,0)$--form. Write $\| \sigma \|_0(x)$ for its comass at $x$. Let $x_0$ be the point at which $\| \sigma \|_0$ attains its maximum. The idea is to construct a simple $(p,0)$--form $\widetilde{\sigma}(x)$ in a neighborhood of $x_0$ such that the $L^2$--norm of $\widetilde{\sigma}$ attains its maximum at $x_0$. See \cite[p. 277]{NiRational} for details. Computing $\partial_v \partial_{\overline{v}} \log \| \widetilde{\sigma} \|^2$ at $x_0$ then yields \begin{eqnarray*}
0 & \geq & \sum_{j=1}^p R_{v \overline{v} j \overline{j}},
\end{eqnarray*}

for all $v \in T_{x_0}^{1,0}M$.

\subsection*{Reminder 3.1}
Fix a point $x \in M$, and let $\Sigma \subset T_x^{1,0}M$ be a $k$--dimensional subspace. Write $\mathbb{S}^{2k-1} \subset \Sigma$ for the unit sphere inside $\Sigma$. Recall that the \textit{$k$--scalar curvature} is defined \begin{eqnarray*}
\text{Scal}_k(x, \Sigma) & : = & k \phantom{.} \dashint_{\mathbb{S}^{2n-1}} \text{Ric}_{\omega}(X, \overline{X}) d \vartheta(X)
\end{eqnarray*}

Using Ni's method of viscosity considerations, we have the following:

\subsection*{Theorem 3.2}
Let $(M^n, \omega)$ be a compact K\"ahler manifold with $\text{Ric}_{\alpha, \beta}^{\perp} >0$ for some $\alpha >0$ and $\beta <0$. If, moreover, $(p+1) \alpha + 2 \beta >0$, then $h^{p,0}=0$ for all $1 \leq p \leq n$. In particular, $M$ is projective.

\begin{proof}
Assume there is a non-zero holomorphic $(p,0)$--form $\sigma \in H_{\overline{\partial}}^{p,0}(M) \simeq H^0(M, \Omega^p(M))$. Let $x_0 \in M$ be the point at which the comass $\| \sigma \|_0$ attains its maximum. From Ni's viscosity considerations, we have (in a fixed unitary frame $e_k$ near $x_0$) \begin{eqnarray}\label{NiViscosity}
\sum_{k=1}^p R_{v \overline{v} k \overline{k}} & \leq & 0
\end{eqnarray}
for any $v \in T_{x_0}^{1,0}M$. Let $\Sigma$ denote the span of $\{ e_1, ..., e_p \}$, and write $\text{Scal}_{p}(x_0, \Sigma)$ for the scalar curvature of $R \vert_{\Sigma}$. The inequality \eqref{NiViscosity} implies that $\text{Scal}_p(x_0, \Sigma) \leq 0$. From $\text{Ric}_{\alpha, \beta}^{\perp}>0$, we observe that \begin{eqnarray*}
0 & < & \dashint_{ X \in \Sigma, | X | =1 } \alpha \text{Ric}_{\omega}(X, \overline{X}) + \beta \text{HSC}_{\omega}(X) d \vartheta(X) \\
&=& \frac{\alpha}{p} \sum_{k=1}^p \text{Ric}_{k \overline{k}} + \frac{2 \beta}{p(p+1)} S_p(x_0, \Sigma) \\
&=& \frac{\alpha}{p} \left( S_p(x_0, \Sigma) + \sum_{\ell =p+1}^n \sum_{k=1}^p R_{\ell \overline{\ell} k \overline{k}} \right) + \frac{2\beta}{p(p+1)} S_p(x_0, \Sigma) \\
& \leq & \frac{\alpha}{p} S_p(x_0, \Sigma) + \frac{2\beta}{p(p+1)} S_p(x_0, \Sigma) \ = \ \frac{(p+1)\alpha + 2\beta}{p(p+1)} S_p(x_0, \Sigma).
\end{eqnarray*}
If $(p+1)\alpha + 2\beta >0$, then $S_p(x_0, \Sigma) >0$, furnishing the desired contradiction.
\end{proof}

Let us raise the following natural question:

\subsection*{Question 3.3}
Let $(M^n, \omega)$ be a compact K\"ahler manifold. Determine all $(\alpha, \beta) \in \mathbb{R}^2$ such that $h^{p,0}=0$ for all $1 \leq p \leq n$. In particular, determine the $(\alpha, \beta) \in \mathbb{R}^2$ such that $M$ is projective and simply connected. \\

To address the simply connectedness problem, let us exhibit the following extension of the diameter estimate in \cite[Proposition 6.2]{ChuLeeTam}:

\subsection*{Proposition 3.4}
Let $(M^n, \omega)$ be a complete K\"ahler manifold with $\text{Ric}_{\alpha, \beta}^{\perp} \geq \lambda$ for some $\alpha > 0 > \beta$, with $3\alpha + 2 \beta>0$ and some $\lambda >0$. Then $M$ is compact with \begin{eqnarray*}
\text{diam}(M,\omega) &\leq & \pi \sqrt{\frac{\alpha(2n-1) + \beta}{\lambda}}.
\end{eqnarray*}

\begin{proof}
Analogous to the argument in the proof of the Bonnet--Meyer theorem, for $p,q \in M$, let $\gamma : [0,\ell] \to M$ be a minimizing geodesic which connects $p$ and $q$. Following \cite{ChuLeeTam}, we will show that \begin{eqnarray*}
\ell & \leq & \pi \sqrt{\frac{\alpha(2n-1) + \beta}{\lambda}}.
\end{eqnarray*}
Let $e_i$ be an orthonormal set of parallel vector fields along $\gamma$ such that $e_{2n-1} = J\dot{\gamma}$ and $e_{2n} = \dot{\gamma}$. Here, $J$ is the underlying complex structure. For each $1 \leq i \leq 2n-1$, set \begin{eqnarray*}
V_i(t) &: =&  \sin \left( \frac{\pi t}{\ell} \right) e_i (t),\\
\phi_i(t,s) &: = & \exp_{\gamma(t)} (s V_i(t)), \\
L_i(s) & : = & \text{length}(\phi_i(\cdot, s)).
\end{eqnarray*}
For each $t$, we have $\phi_i(t,0) = \gamma(t)$, and since $\gamma$ is minimizing, $L_i$ has a minimum at $s=0$. Hence, the second variation of arc length formula yields \begin{eqnarray*}
0 & \leq & \frac{d^2}{ds^2} \Bigg \vert_{s=0} L_i(s) \\
&=& \int_0^{\ell} ( | \nabla V_i |_g^2 - R(V_i, \gamma', \gamma', V_i) ) dt \\
&=& \int_0^{\ell} \left( \left( \frac{\pi}{\ell} \right)^2 \cos^2 \left( \frac{\pi t}{\ell} \right) - \sin^2 \left( \frac{\pi t}{\ell} \right) R(e_i, \gamma', \gamma', e_i) \right) dt.
\end{eqnarray*}
Let $X = \frac{1}{\sqrt{2}}(\gamma' - \sqrt{-1} J \gamma')$. Then $\text{Ric}_{\alpha, \beta}^{\perp}(X) \geq \lambda$ implies \begin{eqnarray*}
\alpha \sum_{i=1}^{2n-1} R(e_i, \gamma', \gamma', e_i) + \beta R(e_{2n-1}, \gamma', \gamma', e_{2n-1}) & \geq & \lambda.
\end{eqnarray*}
Since $3 \alpha + 2 \beta >0$ implies $\alpha(2n-1) + \beta >0$ for $n \geq 2$, we see that \begin{eqnarray*}
0 & \leq & \frac{d^2}{ds^2} \Bigg \vert_{s=0} \left( \alpha \sum_{i=1}^{2n-1} L_i(s) + \beta L_{2n-1}(s) \right) \\
&=& \int_0^{\ell} ( \alpha(2n-1) + \beta) \left( \frac{\pi}{\ell} \right)^2 \cos^2 \left( \frac{\pi t}{\ell} \right) dt - \int_0^{\ell} \sin^2 \left( \frac{\pi t}{\ell} \right)^2 \text{Ric}_{\alpha, \beta}^{\perp}(X) dt \\
& \leq & \int_0^{\ell} \left( (\alpha(2n-1) + \beta) \left( \frac{\pi}{\ell} \right)^2 \cos^2 \left( \frac{\pi t}{\ell} \right) - \lambda \sin^2 \left( \frac{\pi t}{\ell} \right) \right) dt \\
&=& \frac{\ell}{2} \left( (\alpha (2n-1) + \beta) \left( \frac{\pi}{\ell} \right)^2 - \lambda \right).
\end{eqnarray*}
\end{proof}

\subsection*{Remark 3.5}
Note that Chu--Lee--Tam \cite{ChuLeeTam} assume that $\alpha, \beta>0$. Of course, the argument requires only that $\alpha(2n-1) + \beta>0$.

\subsection*{Proposition 3.6}
Let $(M^n, \omega)$ be a compact K\"ahler manifold with $\text{Ric}_{\alpha, \beta}^{\perp}>0$ for some $\alpha > 0 > \beta$ and $\alpha+\beta>0$. Then $M$ is simply connected.

\begin{proof}
Let $\mathcal{O}_M$ denote the structure sheaf of $M$. The Euler characteristic is given by $\chi(\mathcal{O}_M) : = \sum_{k=0}^n (-1)^k h^{k,0}$. For any finite $\nu$--sheeted covering $\widetilde{M} \to M$, the Riemann--Roch--Hirzebruch formula states that \begin{eqnarray*}
\chi(\mathcal{O}_{\widetilde{M}}) &=& \nu \chi(\mathcal{O}_M).
\end{eqnarray*}
If there is a metric on $M$ with $\text{Ric}_{\alpha, \beta}^{\perp}>0$ for some $\alpha >0 > \beta$ and $\alpha+\beta>0$, then $h^{p,0}=0$ for all $1 \leq p \leq n$. From \cite[Theorem 3.2]{NiZhengComparison}, the universal cover $\widetilde{M}$ affords a metric with $\text{Ric}_{\alpha, \beta}^{\perp} > \delta > 0$. In particular, $\widetilde{M}$ is compact and projective with $\chi(\mathcal{O}_{\widetilde{M}})=1$. It follows that $\nu =1$, and thus $\pi_1(M) =0$.
\end{proof}

\subsection*{Conjecture 3.7}
Let $(M^n, \omega)$ be a compact Hermitian manifold. \begin{itemize}
\item[(i)] If $\text{Ric}_{\alpha, \beta}^{(1)} >0$ for some $\alpha > 0$ and $\beta <0$, then $h^{p,0} =0$ for all $1 \leq p \leq n$. 
\item[(ii)] If $\text{Ric}_{\alpha, \beta}^{(2)} >0$ for $\alpha >0$ and $\beta <0$, then $h^{p,0} =0$ for all $1 \leq p \leq n$.
\end{itemize}

\hfill

In this direction, we have the following:

\subsection*{Proposition 3.8}\label{HermitianHodgeEstimates}
Let $(M^n, \omega)$ be a compact Hermitian manifold. \begin{itemize}
\item[(i)] If $\text{Ric}_{\alpha, -\alpha}^{(1)} >0$ for some $\alpha >0$, then $h^{n-1,0} =0$.
\item[(ii)] If $\text{Ric}_{\alpha, \beta}^{(2)} >0$ for some $\alpha >0$, $\beta <0$ with $\alpha + \beta \geq 0$, then $h^{1,0} =0$.
\end{itemize}

\begin{proof}
We observe that Ni's technique of viscosity considerations extends to the Hermitian category. Let $\sigma \in H_{\overline{\partial}}^{n-1,0}(M)$ be a non-trivial holomorphic $(n-1,0)$--form on $M$. Let $p \in M$ be the point at which the comass $\| \sigma \|_0$ achieves its maximum. The Bochner technique, together with the maximum principle implies that \begin{eqnarray}\label{HermBoch1}
\sum_{k=1}^{n-1} R_{v \overline{v} k \overline{k}} & \leq & 0 
\end{eqnarray}
for all $v \in T^{1,0}_p M$. Hence, for any $\alpha >0$, \begin{eqnarray*}
0 \ \geq \ \alpha \sum_{k=1}^{n-1} R_{ n \overline{n} k \overline{k}} &=& \alpha \sum_{k=1}^n R_{ n \overline{n} k \overline{k}} - \alpha R_{n \overline{n} n \overline{n}} \ = \ \text{Ric}_{\alpha, -\alpha}^{(1)},
\end{eqnarray*}
contradicting $\text{Ric}_{\alpha, -\alpha}^{(1)} >0$. This proves (i).

Similarly for (ii), let $\sigma \in H_{\overline{\partial}}^{1,0}(M)$ be a non-trivial holomorphic $(1,0)$--form. The same argument implies (in a fixed unitary frame $e_1, ... , e_n$) \begin{eqnarray}\label{HermBoch2}
R_{v \overline{v}1 \overline{1}} & \leq & 0 
\end{eqnarray}
for any $v \in T_p^{1,0} M$, where $p$ is the point where $| \sigma |^2$ achieves its maximum. Choose $v$ such that $v=e_1$, then \begin{eqnarray*}
\alpha \sum_{k=2}^n R_{k \overline{k} 1 \overline{1}} \ \leq \ 0, \hspace*{1cm} (\alpha + \beta) R_{1 \overline{1} 1 \overline{1}} \ \leq \ 0.
\end{eqnarray*}
Hence, \begin{eqnarray*}
0 & \geq & \alpha \sum_{k=2}^n R_{k \overline{k} 1 \overline{1}} + (\alpha + \beta) R_{1 \overline{1} 1 \overline{1}} \ = \ \text{Ric}_{\alpha, \beta}^{(2)},
\end{eqnarray*}
contradicting $\text{Ric}_{\alpha, \beta}^{(2)}>0$ if $\alpha + \beta \geq 0$. This proves (ii).
\end{proof}

\subsection*{Corollary 3.9}
Let $M^3$ be a compact K\"ahler threefold. If $M$ supports a Hermitian metric with $\text{Ric}_{\alpha, -\alpha}^{(1)} >0$ for some $\alpha >0$, then $M$ is projective. \begin{proof}
From part (i) of \nameref{HermitianHodgeEstimates}, we have $h^{2,0}=0$. It is well-known that a compact K\"ahler manifold satisfying $h^{2,0}=0$ is projective.
\end{proof}

Recall that the Chern scalar curvature of a Hermitian metric $\omega$ is defined by $$\text{Scal}_{\omega} \ : = \ g^{i \overline{j}} g^{k \overline{\ell}} R_{i \overline{j} k \overline{\ell}}.$$ We also have the alterred scalar curvature $\widehat{\text{Scal}}_{\omega} : = g^{i \overline{\ell}} g^{k \overline{j}} R_{i \overline{j} k \overline{\ell}}$. 

\subsection*{Theorem 3.10}
Let $(M^n, \omega)$ be a compact Hermitian manifold with $\text{Ric}_{\alpha, \beta}^{(k)} >0$ for some $1 \leq k \leq 4$, and $\alpha >0 > \beta$. If $\text{Scal}_p(\omega) = \widehat{\text{Scal}}_p(\omega)$ and $(p+1) \alpha + 2 \beta$, then $h^{p,0}=0$ for all $1 \leq p \leq n$.\begin{proof}
Assume there is a non-zero holomorphic $(p,0)$--form $\sigma \in H_{\overline{\partial}}^{p,0}(M) \simeq H^0(M, \Omega^p(M))$. Let $x_0 \in M$ be the point at which the comass $\| \sigma \|_0$ attains its maximum. From Ni's viscosity considerations, we have (in a fixed unitary frame $e_i$ near $x_0$) \begin{eqnarray}\label{Balanced}
\sum_{i=1}^p R_{v \overline{v} i \overline{i}} & \leq & 0 
\end{eqnarray}
for any $v \in T_{x_0}^{1,0}M$. Let $\Sigma : = \text{span} \{ e_1, ..., e_p \}$, and write $\text{Scal}_p(x_0, \Sigma)$ for the Chern scalar curvature of $R \vert_{\Sigma}$. Similarly, write $\widehat{\text{Scal}}_p(x_0, \Sigma)$ for the alterred scalar curvature of $R \vert_{\Sigma}$. If $\omega$ is balanced, then $\text{Scal}_p(x_0, \Sigma) = \widehat{\text{Scal}}_p(x_0, \Sigma)$, and \eqref{Balanced} implies $\text{Scal}_p(x_0, \Sigma) \leq 0$. Assume $\text{Ric}_{\alpha,\beta}^{(k)} >0$ for some $1 \leq k \leq 3$, then \begin{eqnarray*}
0 & < & \dashint_{\mathbb{S}^{2p-1}} \alpha \text{Ric}_{\omega}^{(k)}(X, \overline{X}) + \beta \text{HSC}_{\omega}(X) d\vartheta \\
&=& \frac{\alpha}{p}\sum_{i=1}^p \text{Ric}_{i \overline{i}}^{(k)} + \frac{2\beta}{p(p+1)}(\text{Scal}_p(x_0, \Sigma) + \widehat{\text{Scal}}_p(x_0, \Sigma))\\
&=& \frac{\alpha}{p} \text{Scal}_p(x_0, \Sigma) + \frac{2\beta}{p(p+1)} (\text{Scal}_p(x_0, \Sigma) + \widehat{\text{Scal}}_p(x_0, \Sigma) ) \\
&=& \frac{(\alpha(p+1) + \beta) \text{Scal}_p(x_0, \Sigma) + \beta \widehat{\text{Scal}}_p(x_0, \Sigma)}{2p(p+1)} \\
&=& \frac{\alpha(p+1)+2\beta}{2p(p+1)} \text{Scal}_p(x_0, \Sigma),
\end{eqnarray*}
where the last equality makes use of the balanced condition. If $\alpha(p+1)+2\beta>0$, we have the desired contradiction.
\end{proof}

From \cite{LiuYangRicci}, we have the pointwise equality \begin{eqnarray*}
\text{Scal}_{\omega} & = & \widehat{\text{Scal}}_{\omega}  + \langle \overline{\partial} \overline{\partial}^{\ast} \omega, \omega \rangle.
\end{eqnarray*}

In particular, if $\omega$ is balanced, then $\text{Scal}_{\omega} = \widehat{\text{Scal}}_{\omega}$:

\subsection*{Corollary 3.11}
Let $(M^n, \omega)$ be a compact Hermitian manifold with a balanced metric of $\text{Ric}_{\alpha, \beta}^{(k)} >0$ for some $1 \leq k \leq 4$, and $\alpha > 0 > \beta$. If $\alpha(n+1) + 2 \beta >0$, then $h^{n,0}=0$.\\

Of course, when $n=2$, i.e., on compact complex surfaces, the balanced condition is equivalent to the K\"ahler condition. 

\subsection*{Remark 3.12}
Petersen--Wink \cite{PetersenWinkHodge} have established estimates on the Hodge numbers of K\"ahler manifolds in terms of the eigenvalues of the K\"ahler curvature operator. In particular, they show that if $\mathfrak{R} : \mathfrak{u}(n) \to \mathfrak{u}(n)$ denotes the K\"ahler curvature operator with eigenvalues $\lambda_1 \leq \lambda_2 \leq \cdots \leq \lambda_{n^2}$,  then $h^{2,0}=0$ if $$\lambda_1 + \cdots + \lambda_{n-1} > 0.$$ That is, $M$ is projective if the metric is $(n-1)$--positive. Let us remark that is quite strong: for $n=3$,  this reduces to the $2$--positivity of the K\"ahler curvature operator, which is known to be equivalent to the positivity of the orthogonal bisectional curvature $\text{HBC}_{\omega}^{\perp}$. In particular, such metrics are all biholomorphically isometric to $(\mathbb{P}^3, \omega_{\text{FS}})$.

If $(M, \omega)$ is a compact K\"ahler surface, the orthogonal Ricci curvature $\text{Ric}_{\omega}^{\perp}$ is equivalent to the orthogonal bisectional curvature $\text{HBC}_{\omega}^{\perp}$. From the work of Gu--Zhang \cite{GuZhang}, any compact simply connected K\"ahler surface with $\text{Ric}_{\omega}^{\perp} >0$, therefore, deforms under the K\"ahler--Ricci flow to a metric with $\text{HBC}_{\omega} >0$. It follows from the solution of the Frankel conjecture that $M$ is biholomorphic to $\mathbb{P}^2$.

\subsection*{Proposition 3.13}
Let $(M^n, \omega)$ be a compact K\"ahler manifold with $\text{Ric}_{\alpha, \beta}^{\perp} <0$ for some $\alpha >0$ and $\beta <0$. If $\alpha + \beta \geq 0$, then $M$ does admit any non-trivial holomorphic vector fields. \begin{proof}
Let $X$ be a non-trivial holomorphic vector field. The Bochner formula gives \begin{eqnarray*}
\left \langle \sqrt{-1} \partial \overline{\partial} | X |^2, \frac{1}{\sqrt{-1}} v \wedge \overline{w} \right \rangle &=& \langle \nabla_v X, \overline{\nabla_w X} \rangle - R_{v \overline{w} X \overline{X}}.
\end{eqnarray*}
At the point $p \in M$ where $| X |^2$ attains its non-zero maximum, we have \begin{eqnarray}\label{VecBoch}
R_{v \overline{v} X \overline{X}} & \geq & 0, 
\end{eqnarray}
for all $v \in T_p^{1,0}M$. Let $\{ e_1, ..., e_n \}$ be a local unitary frame near $p$, such that $X \vert_p = e_1 \vert_p$. From \eqref{VecBoch}, we have \begin{eqnarray}\label{VecBoch2}
\alpha \left( \sum_{k=2}^n R_{k \overline{k} 1 \overline{1}} \right) \ \geq \ 0, \hspace*{1cm} (\alpha + \beta) \text{HSC}_{\omega}(X) \ \geq \ 0.
\end{eqnarray}
Since $\text{Ric}_{\alpha, \beta}^{\perp} <0$, however, \begin{eqnarray*}
\alpha \left( \sum_{k=2}^n R_{k \overline{k} 1 \overline{1}} \right) + (\alpha + \beta) \text{HSC}_{\omega}(X) &< & 0,
\end{eqnarray*}
violating \eqref{VecBoch2}.
\end{proof}

Immediate from the argument in the K\"ahler category, is the following: 

\subsection*{Proposition 3.14}
Let $(M^n, \omega)$ be a compact Hermitian manifold. Suppose $\text{Ric}_{\alpha, \beta}^{(2)} <0$ for $\alpha >0$, $\beta <0$, with $\alpha + \beta \geq 0$. Then $M$ does not admit any non-trivial holomorphic vector fields.\\

Let us close this section by extending an old result of Cheung \cite{Cheung}: 

\subsection*{Theorem 3.15}
Let $(M^2, \omega)$ be a K\"ahler--Einstein surface with $\text{Ric}_{\omega} = \lambda \omega$. The metric $\omega$ has negative holomorphic sectional curvature if and only if \begin{eqnarray*}
\text{Ric}_{2,-1}^{\perp} <0 \hspace*{1cm} \text{and} \hspace*{1cm} | R_{1 \overline{2} 1 \overline{2}} |^2 \ < \ | \text{Ric}_{2, -1}^{\perp} |^2.
\end{eqnarray*}

\begin{proof}
Fix a point $p \in M$, and assume the holomorphic sectional curvature achieves a minimum in the direction $e_1$. Taking partial derivatives of $\sum_{i,j,k,\ell} R_{i \overline{j} k \overline{\ell}} v_i \overline{v}_j v_k \overline{v}_{\ell}$, we see that $R_{1 \overline{1} 1 \overline{2}} = R_{1 \overline{1} 2 \overline{1}} =0$ at $p$. Since the metric is K\"ahler--Einstein, we further deduce that the following components of the curvature vanish: \begin{eqnarray*}
R_{1 \overline{1} 1 \overline{2}} \ = \ R_{1 \overline{1} 2 \overline{1}} \ = \ R_{1 \overline{2} 1 \overline{1}} \ = \ R_{1 \overline{2} 2 \overline{2}} \ = \ R_{2 \overline{1} 1 \overline{1}} \ = \ R_{2 \overline{1} 2 \overline{2}} \ = \ R_{2 \overline{2} 1 \overline{2}} \ = \ R_{2 \overline{2} 2 \overline{1}} \ = \ 0.
\end{eqnarray*}
The holomorphic sectional curvature at $p$, in the unit direction $(v_1, v_2)$ is given by \begin{eqnarray*}
\sum_{i,j,k,\ell} R_{i \overline{j} k \overline{\ell}} v_i \overline{v}_j v_k \overline{v}_{\ell} &=& R_{1 \overline{1} 1 \overline{1}} (v_1 \overline{v}_1)^2 + 4 R_{1 \overline{1} 2 \overline{2}} v_1 \overline{v}_1 v_2 \overline{v}_2 + R_{2 \overline{1} 2 \overline{1}} (v_2 \overline{v}_1)^2 \\ 
&& \hspace*{3cm} + R_{1 \overline{2} 1 \overline{2}} (v_1 \overline{v}_2)^2 + R_{2 \overline{2} 2 \overline{2}}(v_2 \overline{v}_2)^2 \\ \\
&=& R_{1 \overline{1} 1 \overline{1}} + 2(2R_{1 \overline{1} 2 \overline{2}} - R_{1 \overline{1} 1 \overline{1}}) | v_1 \overline{v}_2 |^2 +2 \text{Re} \left( R_{1 \overline{2} 1 \overline{2}} (v_1 \overline{v}_2)^2 \right).
\end{eqnarray*}
Write $R_{1 \overline{2} 1 \overline{2}} = | R_{1 \overline{2} 1 \overline{2}} | e^{i \vartheta_1}$ and consider the direction $v_1 = \frac{1}{\sqrt{2}}$, $v_2 = \frac{\sqrt{-1}}{\sqrt{2}} e^{i \vartheta_1/2}$. The holomorphic sectional curvature in this direction is therefore \begin{eqnarray*}
&& R_{1 \overline{1} 1 \overline{1}} + \frac{1}{2} (2R_{1 \overline{1} 2 \overline{2}} - R_{1 \overline{1} 1 \overline{1}}) +2 \text{Re} \left( -\frac{1}{4} | R_{1 \overline{2} 1 \overline{2}} | e^{i \vartheta_1} e^{-i \vartheta_1} \right) \\
&& \hspace*{6cm} =  R_{1 \overline{1} 1 \overline{1}} + \frac{1}{2}(2R_{1 \overline{1} 2 \overline{2}} - R_{1 \overline{1} 1 \overline{1}} - | R_{1 \overline{2} 1 \overline{2}} | ).
\end{eqnarray*}
Since $e_1$ minimizes the holomorphic sectional curvature, we have \begin{eqnarray*}
2 R_{1 \overline{1} 2 \overline{2}} - R_{1 \overline{1} 1 \overline{1}} & \geq & | R_{1 \overline{2} 1 \overline{2}} |,
\end{eqnarray*}
which implies $2 R_{1 \overline{1} 2 \overline{2}} \geq R_{1 \overline{1} 1 \overline{1}}$, i.e., $2 \lambda \geq 3 R_{1 \overline{1} 1 \overline{1}}$. 

Extending the above calculation, \begin{eqnarray*}
\sum_{i,j,k,\ell} R_{i \overline{j} k \overline{\ell}} v_i \overline{v}_j v_k \overline{v}_{\ell} & \leq & R_{1 \overline{1} 1 \overline{1}} + 2(2 R_{1 \overline{1} 2 \overline{2}} - R_{1 \overline{1} 1 \overline{1}}) | v_1 \overline{v}_2 |^2 + 2 | R_{1 \overline{2} 1 \overline{2}} | | v_1 \overline{v}_2 |^2 \\
& = & R_{1 \overline{1} 1 \overline{1}} + 2(2R_{1 \overline{1} 2 \overline{2}} - R_{1 \overline{1} 1 \overline{1}} + | R_{1 \overline{2} 1 \overline{2}} |) | v_1 \overline{v}_2 |^2 \\
& \leq & R_{1 \overline{1} 1 \overline{1}} + \frac{1}{2} (2R_{1 \overline{1}2 \overline{2}} - R_{1 \overline{1} 1 \overline{1}} + | R_{1 \overline{2} 1 \overline{2}} |) ( | v_1 |^2 + | v_1 |^2)^2 \\
&=& R_{1 \overline{1} 1 \overline{1}} + \frac{1}{2}(2 R_{1 \overline{1} 2 \overline{2}} - R_{1 \overline{1} 1 \overline{1}} + | R_{1 \overline{2} 1 \overline{2}} |).
\end{eqnarray*}
A variation argument similar to the one above shows that the upper bound is achieved when $v_1 = \frac{1}{\sqrt{2}}$ and $v_2 = \frac{1}{\sqrt{2}} e^{i \vartheta_1/2}$. Hence, the holomorphic sectional curvature is maximized at $p$ with value $$R_{1 \overline{1} 1 \overline{1}} + \frac{1}{2}(2 R_{1\overline{1} 2 \overline{2}} - R_{1 \overline{1} 1 \overline{1}} + | R_{1 \overline{2} 1 \overline{2}} |).$$ Since this equantity is negative if and only if $$\lambda - \frac{1}{2}R_{1 \overline{1} 1 \overline{1}} + \frac{1}{2} | R_{1 \overline{2} 1 \overline{2}} | < 0$$ if and only if $$2 \lambda - R_{1 \overline{1} 1 \overline{1}} < 0 \hspace*{1cm} \text{and} \ | R_{1 \overline{2} 1 \overline{2}} |^2 < | 2 \lambda - R_{1 \overline{1} 1 \overline{1}} |^2,$$ this completes the proof.
\end{proof}

\section{Examples}

\subsection*{Example 4.1. The Iwasawa Threefold}
Let $X = G / \Gamma $ denote the Iwasawa threefold given by the quotient of \begin{eqnarray*}
G & : =& \left \{ \begin{pmatrix}
1 & z_1 & z_3 \\
0 & 1 & z_2 \\
0 & 0 & 1
\end{pmatrix} : (z_1, z_2, z_3) \in \mathbb{C}^3 \right \}
\end{eqnarray*}

by the discrete group \begin{eqnarray*}
\Gamma & : = & \left \{ \begin{pmatrix}
1 & z_1 & z_3 \\
0 & 1 & z_2 \\
0 & 0 & 1
\end{pmatrix} : z_1, z_2, z_3 \in \mathbb{Z} + \sqrt{-1} \mathbb{Z} \right \}.
\end{eqnarray*}

It is well-known that $X$ is non-K\"ahler, but supports a balanced metric. Indeed, the projection map $f : X \to \mathbb{Z}[\sqrt{-1}]$ given by $$f : \begin{pmatrix}
1 & z_1 & z_3 \\
0 & 1 & z_2 \\
0 & 0 & 1
\end{pmatrix} \ \mapsto \ z_1$$ is a surjective holomorphic map with K\"ahler fibers. The map $$\sigma : 
z_1 \ \mapsto \  \begin{pmatrix}
1 & z_1 & z_3 \\
0 & 1 & z_2 \\
0 & 0 & 1
\end{pmatrix}$$ defines a holomorphic section of $f$. By \cite[Theorem 5.5]{Michelson}, $X$ admits a balanced metric. The Hodge numbers of $X$ are detailed in \cite[p. 49]{DanielleAngela}. In particular, \begin{eqnarray*}
h^{1,0} = 3, \hspace{0.5cm} h^{2,0} = 3, \hspace{0.5cm} h^{3,0} = 1.
\end{eqnarray*}

We, therefore, have the following:

\subsection*{Corollary 4.2}
Let $X$ be the Iwasawa threefold. There is no balanced Hermitian metric on $X$ with $\text{Ric}_{\alpha, \beta}^{(k)} >0$ for $\alpha > 0 > \beta$ and $3 \alpha + 2 \beta >0$.

\subsection*{Example 4.3.  $\text{U}(n)$--invariant K\"ahler metrics on $\mathbb{C}^n$}

In the standard coordinates on $\mathbb{C}^{n \geq 3}$, a $\text{U}(n)$--invariant K\"ahler metric is given by \begin{eqnarray*}
g_{i \overline{j}} &=& f(r) \delta_{ij} + f'(r) \overline{z}_i z_j,
\end{eqnarray*}

where $r = \sum_{k=1}^n | z_k |^2$, and the function $f$ is smooth on $[0,\infty)$. Set $h = (rf)'$.\\

We first recall the following lemma of Wu--Zheng \cite{WuZheng}:

\subsection*{Lemma 4.4}
The $\text{U}(n)$--invariant metric $g$ defined above is a complete K\"ahler metric if and only if $f >0$, $h>0$, and \begin{eqnarray*}
\int_0^{\infty} \sqrt{\frac{h}{r}} dr &=& \infty.
\end{eqnarray*}

For $h>0$, the function $\xi = - rh'/h$ is smooth on $[0,\infty)$, with $\xi(0)=0$. The components of the curvature tensor of a $\text{U}(n)$--invariant K\"ahler metric in the unitary frame $$e_1 : = \frac{1}{\sqrt{h}} \partial_{z_1}, \hspace{1cm} e_k : = \frac{1}{\sqrt{f}} \partial_{z_k}, \hspace{0.5cm} k \geq 2,$$ at the point $p = (z_1, 0, ..., 0)$ are given by \begin{eqnarray*}
A & : = & R_{1 \overline{1} 1 \overline{1}} \ = \ \xi'/h \\
B & : = & R_{1 \overline{1} i \overline{i}} \ = \ \frac{1}{(rf)^2} \left[ rh - (1-\xi) \int_0^r h ds \right], \hspace*{1cm} i \geq 2, \\
C & : = & R_{i \overline{i} i \overline{i}} \ = \ 2 R_{i \overline{i} j \overline{j}} \ = \ \frac{2}{r^2 f^2} \left( \int_0^r h ds - rh \right), \hspace*{1cm} 2 \leq i \neq j.
\end{eqnarray*}

All other components (not given by symmetries of the above) vanish. By the unitary-invariance, it suffices to calculate the curvature at the point $p$.

In the above notation (c.f., \cite[p. 6]{HuangTam}), we have \begin{eqnarray*}
R_{1 \overline{1}} &=& A + (n-1) B, \\
R_{1 \overline{1}1 \overline{1}} &=& R_{1 \overline{1}} - (n-1)B \\
R_{i \overline{i}} &=& B + C + \frac{(n-2)}{2} C \ = \ B + \frac{n}{2} C, \hspace*{1cm} i \geq 2, \\
R_{i \overline{i} i \overline{i}} &=& R_{i \overline{i}} - B - \frac{(n-2)}{2}C, \hspace*{1cm} i \geq 2.
\end{eqnarray*}

\subsection*{Proposition 4.5}
Let $\omega$ be the complete $\text{U}(n)$--invariant metric above. In the above notation, \begin{eqnarray*}
\text{Ric}_{\alpha,\beta}^{\perp}(e_1) &=& \alpha(n-1) B + \beta A, 
\end{eqnarray*}
and for each $i \geq 2$, 
\begin{eqnarray*}
\text{Ric}_{\alpha,\beta}^{\perp}(e_i) &=&\left( \frac{\alpha n}{2} + \beta \right) C + \alpha B.
\end{eqnarray*}

\subsection*{Example 4.6. Hopf Manifolds}
Let $X = \mathbb{S}^{2n-1} \times \mathbb{S}^1$ be the standard Hopf manifold of dimension $n \geq 2$. On $X$ there is a natural metric (inheritted from the cyclic group action of $z \mapsto \frac{1}{2} z$ on $\mathbb{C}^n - \{ 0 \}$): \begin{eqnarray*}
\omega_0 &=& \sqrt{-1} g_{i \overline{j}} dz_i \wedge d\overline{z}_j \ = \ \sqrt{-1} \frac{4\delta_{ij}}{| z |^2} dz_i \wedge d\overline{z}_j.
\end{eqnarray*}

From \cite[$\S 6.1$]{LiuYangRicci}, the components of the Chern curvature tensor read \begin{eqnarray*}
R_{i \overline{j} k \overline{\ell}} &=& g^{p \overline{q}} \frac{\partial g_{p \overline{\ell}}}{\partial z_i} \frac{\partial g_{k \overline{q}}}{\partial \overline{z}_j} - \frac{\partial^2 g_{k \overline{\ell}}}{\partial z_i \partial \overline{z}_j} \ = \ \frac{4\delta_{k\ell}(\delta_{ij} | z |^2 - z_j \overline{z}_i)}{| z |^6}.
\end{eqnarray*}

Let $v \in T^{1,0}X$ be a $(1,0)$--tangent vector of unit length. Then \begin{eqnarray*}
\text{HSC}_{\omega}(v) & = & \sum_{i,j,k,\ell=1}^n \frac{4\delta_{k\ell}(\delta_{ij} | z |^2 - z_j \overline{z}_i)}{| z |^6} v_i \overline{v}_j v_k \overline{v}_{\ell} \\ 
&=&  \frac{4}{| z |^6} \sum_{k=1}^n  | v_k |^2 \sum_{i \neq j =1}^n \left( |v_i |^2 | z |^2 - z_j \overline{z}_i v_i \overline{v}_j \right).
\end{eqnarray*}

Moreover, we have that \begin{eqnarray*}
\text{Ric}_{\omega}^{(1)}(v) \ = \ n \sqrt{-1} \partial \overline{\partial} \log | z |^2 
\end{eqnarray*}

and \begin{eqnarray*}
\text{Ric}_{\omega}^{(2)}(v) \ = \ \frac{n-1}{4} \left( \frac{4}{| z |^2} \sum_{i=1}^n | v_i |^2 \right) \ = \ \frac{n-1}{| z |^2}.
\end{eqnarray*}

Hence, \begin{eqnarray*}
\text{Ric}_{\alpha, \beta}^{(2)} &=& \frac{\alpha(n-1)}{| z |^2} + \frac{4\beta}{| z |^6} \sum_{k=1}^n \sum_{i \neq j=1}^n ( |v_i |^2 | z |^2 - z_j \overline{z}_i v_i \overline{v}_j )
\end{eqnarray*}

\subsection*{Acknowledgements}
The first author would like to thank his advisors Ben Andrews and Gang Tian for their unwavering support and encouragement; and Gang Tian for corrections on the earlier version of this manuscript. The second author is grateful to Professor Fangyang Zheng for constant encouragement and support.

\hfill

\scshape{Kyle Broder. Mathematical Sciences Institute, Australian National University, Acton, ACT 2601, Australia}

\scshape{BICMR, Peking University, Beijing, 100871, People's republic of china}

\textit{E-mail address}: \texttt{Kyle.Broder@anu.edu.au} 

\hfill

\scshape{Kai Tang. College of Mathematics and Computer Science, Zhejiang Normal University, Jinhua, Zhejiang, 321004, China}

\textit{Email address:} \texttt{kaitang001@zjnu.edu.cn}

\end{document}